# The Multi-set Allocation Occupancy function and inequality (MAO function and MAO inequality): the foundation of Generalized hypergeometric distribution theory

Xing-gang Mao[a*] and Xiao-yan Xue[b]

[a]Department of Neurosurgery, Xijing Hospital, the Fourth Military Medical University, Xi'an, No. 17 Changle West Road, Xi'an, Shaanxi Province, China;

[b]Department of Pharmacology, School of Pharmacy, the Fourth Military Medical University, No. 17 Changle West Road, Xi'an, Shaanxi Province, China

*Corresponding authors. E-mail addresses: xgmao@fmmu.edu.cn (X.G. Mao)

Short title: MAO function, inequality and distribution

**Abstract**

In our previous work, we studied the Generalized hypergeometric distribution (GHGD), which we refer to as multi-set allocation occupancy (MAO) distribution. We derived and proved formulas for the mathematical expectation and variance of the GHGD for **any number of subsets $T$ and any overlap count $t$ (where $1 \leq t \leq T$)**. In addition, we give an Asymptotic property of GHGD. However, the formulas as reported in our previous work exhibited an unsatisfactory complex form, and the formulas of higher moments of the GHGD are not derived. After further deeply study, we established a novel, unreported function, which describes the higher moments of the MAO distributions with an astonishing unified, elegant, and concise formula. The core definitions include: **multi-set allocation occupancy (MAO) function** $g(A_1, A_2, \ldots, A_r) \prod_{i=1}^{T} [(m_i)_{k_i} \cdot (n - m_i)_{r-k_i}]$ **and MAO norm:** $\|(p_1, \ldots, p_r)\|_T = \frac{\sum_{\substack{A_1, \ldots, A_r \subseteq [T] \\ |A_1|=p_1, \ldots, |A_r|=p_r}} g(A_1, A_2, \ldots, A_r)}{((n)_r)^{T-1}}$, where $p_i$ is the size of the subset $A_i$ from the index vector $[T] = \{1, 2, \ldots, T\}$, $m_i < n$, and $(x)_r$ is the falling factorial. With this new definition, the intricate moment relations collapse into an unified form: any raw moment ($v$th) of $p(x_{=t})$ and $p(x_{\geq t})$ can be calculated by the formulas:

$$E(x_{=t}^v) = \sum_{1 \leq i \leq v} s_{v,i} \|t^i\|$$

And

$$E(x_{\geq t}^v) = \sum_{1 \leq i \leq v} s_{v,i} \|[t, T]^i\|$$

Where $s_{v,i}$ are Stirling numbers of the second kind and $[t, T] = \{t, t+1, \ldots, T\}$. In addition, interestingly, based on the MAO norm definition, we get the MAO inequality on the condition that all $p_i$ satisfy the proximity condition: $max(p_i) - min(p_i) \leq 1$, a relatively conservative condition):

$$\prod_{1 \leq i \leq r} \|(p_i)\|_T \geq \|(p_1, \ldots, p_r)\|_T$$

As a direct corollary, we can get the asymptotic property of MAO distribution that $E(x) > Var(x)$ and $E(x) - Var(x) = o(E(x))$ when $E(x) \to 0$.

1. Background

*General hypergeometric distribution* (**GHGD**), or **multi-set allocation occupancy (MAO) distribution**: from a finite space $\Omega$ containing $n$ elements, randomly select totally $T$ subsets $S_i$ (each contains $m_i$ elements, $1 \leq i \leq T$; here $m_i = s_i = |S_i|$, and $\boldsymbol{m}$ is the count vector $\boldsymbol{m} = \{m_1, m_2, \ldots m_T\}$), what is the probability that exactly $x$ elements are overlapped exactly $t$ times or at least $t$ times? Here, the random variables are $x_{=t}$ and $x_{\geq t}$. We have provide the mathematical expectation and variance for the following distributions in previous manuscripts:

$$p(x_{=t} = x \mid n, \{m_1, m_2 \ldots m_T\}) \quad \text{or} \quad p(x_{=t} = k \mid n, \boldsymbol{m}) \quad \text{or} \quad p(x_{=t})$$

and

$$p(x_{\geq t} = k \mid n, \{m_1, m_2 \ldots m_T\}) \quad \text{or} \quad p(x_{\geq t} = k \mid n, \boldsymbol{m}) \quad \text{or} \quad p(x_{\geq t})$$

We have comprehensively resolved the mathematical expectation and variance of the GHGD for any $t$ values. However, the formulas as reported in our previous work exhibited an unsatisfactory complex form [1,2]. After deeply study the formulas and higher moments, we established a novel, unreported function and inequality, which provides core definitions to describe the higher moments of the above distributions with an astonishing unified, elegant, and concise formula.

First, we define the $k$th descending factorial of integer $m_i$, and $(m_i)_0 = 1$:

$$(m_i)_k = \prod_{j=0}^{k-1} (m_i - j)$$

Second, define $A$, $B$, …, or $A_1$, $A_2$…, as subsets of index vector $[T]$ ($[T] = \{1, 2, \ldots, T\}$). Then $m_A = \{m_i \mid i \in A\}$ is a subset of vector $\boldsymbol{m}$. If the number of elements in $A$ has to be emphasized it can be denoted as $A_{|A|=t}$, indicating the cardinality of $A$ is $t$.

## 2. Definition of Multi-set Allocation Occupancy function (MAO function)

Let $[T] = \{1, 2, \dots, T\}$ be a set of indices. For each index $i \in [T]$, let $m_i$ be a positive integer with $1 \leq m_i \leq n$, where $n$ is a positive integer.

For any subset $A \subseteq [T]$, define the function $g(A)$ as:

$$g(A) = \prod_{i \in A} m_i \cdot \prod_{i \notin A} (n - m_i)$$

This function represents the product over all indices, where each index $i$ contributes $m_i$ if $i \in A$, and $n - m_i$ if $i \notin A$.

### 2.1. Define the falling factorial $(x)_r$ for $r \geq 1$, and $(x)_0 = 1$:

$$(x)_r = x(x-1) \cdots (x - r + 1)$$

### 2.2. Define the Allocation Occupancy function (MAO function) $g(A_1, A_2, \dots, A_r)$:

For a collection of $k$ subsets $A_1, A_2, \dots, A_r \subseteq [T]$, define the function $g(A_1, A_2, \dots, A_r)$ as:

$$g(A_1, A_2, \dots, A_r) = \prod_{i=1}^{T} [(m_i)_{k_i} \cdot (n - m_i)_{r - k_i}]$$

where $k_i$ is the number of times index $i$ appears in the subsets $A_1, A_2, \dots, A_r$, i.e., $k_i = \sum_{j=1}^{r} 1_{A_j}$, and $(x)_r$ denotes the falling factorial.

### 2.3. Special forms of the MAO function.

If all $A_i$ is the same subset $A$, we have:

$$g(A^r) = \prod_{i \in A} (m_i)_r \prod_{i \notin A} (n - m_i)_r$$

When all $A = [T]$, we have:

$$g([T]^r) = \prod_{i \in [T]} (m_i)_r$$

Especially,

$$g([T]) = \prod_{i \in [T]} m_i$$

When all $A = \emptyset$, we have:

$$g(\emptyset^r) = \prod_{i \in [T]} (n - m_i)_r$$

Especially,

$$g(\emptyset) = \prod_{i \in [T]} (n - m_i)$$

When there is only one sunset $A$, we can define this special MAO function as:

$$f(A) = g(A) = \prod_{i \in A} m_i \prod_{i \notin A} (n - m_i)$$

Therefore, the operation function of $g(A_1, A_2, \ldots, A_r)$ is a generalized function for all situations.

## 3. Definition of Multi-set Allocation Occupancy norm (MAO norm)

3.1. Define the transversal sum of $g(A_1, A_2, \ldots, A_r)$ as:

$$G_T(p_1, \ldots, p_r) = \sum_{\substack{A_1, \ldots, A_r \subseteq [T] \\ |A_1| = p_1, \ldots, |A_r| = p_r}} g(A_1, A_2, \ldots, A_r)$$

3.2. Define the most general type of $g(A_1, A_2, \ldots, A_r)$ as:

$$G_T(B_1, \ldots, B_r) = \sum_{\substack{A_1, \ldots, A_r \subseteq [T] \\ |A_1| \in B_1, \ldots, |A_r| \in B_r \\ B_1, \ldots, B_r \subseteq [T]}} g(A_1, A_2, \ldots, A_r)$$

Therefore, $G_T(p_1, \ldots, p_r)$ can be viewed as special $G_T(B_1, \ldots, B_r)$, where all $B_i$ has only one element $p_i$.

3.3. Define the norm of $p_1, p_2, \ldots, p_r$ as

$$\|(p_1, \ldots, p_r)\|_T = N_T(p_1, \ldots, p_r) = \frac{G_T(p_1, \ldots, p_r)}{((n)_r)^{T-1}}$$

$$= \frac{\sum_{\substack{A_1, \ldots, A_r \subseteq [T] \\ |A_1| = p_1, \ldots, |A_r| = p_r}} g(A_1, A_2, \ldots, A_r)}{((n)_r)^{T-1}}$$

3.4. Define the most generalized definition as:

$$\|(B_1, \ldots, B_r)\|_T = N_T(B_1, \ldots, B_r) = \frac{G_T(B_1, \ldots, B_r)}{((n)_r)^{T-1}}$$

3.5. Define the vector form of $G_T$ and MAO norm:

$$G_T(\mathbf{p})$$

$$\|\mathbf{p}\|_T = N_T(\mathbf{p})$$

Where $\mathbf{p} = \{P_1, P_2, \ldots, P_r\}$, and $P_i$ are elements or subsets of $[T]$ ($P_i \in [T]$ or $P_i \subseteq [T]$).

**Theorem 1 (MAO theory for GHGD)**

Any raw moment of $p(x_{=t})$ and $p(x_{\geq t})$ can be calculated by the unified formulas:

$$E(x_{=t}^v) = \sum_{1 \leq i \leq v} s_{v,i} \|t^i\|$$

And

$$E(x_{\geq t}^v) = \sum_{1 \leq i \leq v} s_{v,i} \|[t,T]^i\|$$

Where $s_{v,i}$ are Stirling numbers of the second kind and $[t,T] = \{t, t+1, \ldots, T\}$. The first several commonly used $s_{v,i}$ were listed in Table S2 for convenience.

4. Statement of the MAO Inequality

**Theorem 2 (MAO Inequality).** Let $p_1, p_2, \ldots, p_r$ be non-negative integers with $0 \leq p_j \leq T$ for all $j$. Then the following inequality holds on the condition that all $p_i$ satisfy the proximity condition: $max(p_i) - min(p_i) \leq 1$, which is a relatively conservative condition (a more precise condition is $max(p_i) - min(p_i) \leq max(1, r-2)$):

$$\sum_{\substack{A_1,\ldots,A_r \subseteq [T] \\ |A_1|=p_1,\ldots,|A_r|=p_r}} \frac{\prod_{j=1}^{r} g(A_j)}{n^{r(T-1)}} \geq \sum_{\substack{A_1,\ldots,A_r \subseteq [T] \\ |A_1|=p_1,\ldots,|A_r|=p_r}} \frac{g(A_1, A_2, \ldots, A_r)}{((n)_r)^{T-1}}.$$

Given that:

$$\sum_{\substack{A_1,\ldots,A_r \subseteq [T] \\ |A_1|=p_1,\ldots,|A_r|=p_r}} \frac{\prod_{j=1}^{r} g(A_j)}{n^{r(T-1)}} = \prod_{i=1}^{r} \sum_{\substack{A_i \subseteq [T] \\ |A_i|=p_i}} \frac{g(A_i)}{n^{(T-1)}}$$

The MAO inequality can be abbreviated as the norm form:

$$\prod_{1 \leq i \leq r} \|(p_i)\|_T \geq \|(p_1, \ldots, p_r)\|_T$$

Or

$$\prod_{i=1}^{r} N_T(p_i) \geq N_T(p_1, \ldots, p_r)$$

Or expressed as the vector form ($\mathbf{p} = \{p_1, \ldots, p_r\}$):

$$\prod_{1 \leq i \leq r} \|(p_i)\|_T \geq \|\mathbf{p}\|_T$$

Strikingly and interestingly, the norm form of the MAO inequality shares exactly the same form of Hölder-Type Inequality ( $\| f_1 f_2 \cdots f_r \|_1 \leq \| f_1 \|_{p_1} \| f_2 \|_{p_2} \cdots \| f_r \|_{p_r}$, $\sum_{i=1}^{r} \frac{1}{p_i} = 1$), implying a deep link between the normed space and set-valued space. The MAO inequality indicated that the norm of the product of functions is dominated by the product of their individual norms.

This analogy reveals a deep mathematical principle: **in both continuous function spaces and discrete set systems, the collective behavior is controlled by the product of individual behaviors**. This underscores a universal theme in mathematics where independence or separability leads to such inequalities.

The MAO inequality can thus be seen as a combinatorial analogue of Hölder's inequality, extending the concept of norm domination to discrete set-values spaces. This connection suggests that MAO function and inequality may serve as a foundational tool in combinatorial analysis, with potential applications in probability, optimization, and theoretical computer science.

## 5. Proof of the MAO Inequality

We will give the proof of the most common condition, that all $p_i = p$. The detailed proof is provided in the Appendix.

## 6. Application in GHGD study

### 6.1. Corollary: Asymptotic property of GHGD: $E(x) > Var(x)$ and $E(x) - Var(x) = o(E(x))$ when $E(x) \to 0$.

Let:

$$\Delta_{EV} = E(x) - Var(x)$$

where $x$ can be either $x_{=t}$ or $x_{\geq t}$ for any $t$ value.

According to the MAO theorem, for $x = x_{=t}$ we have:

$$E(x) = \|t\|_T$$

$$Var(x) = E(x^2) - (E(x))^2 = \|t\|_T + \|t^2\|_T - (\|t\|_T)^2$$

$$\Delta_{EV} = E(x) - Var(x) = (\|t\|_T)^2 - \|t^2\|_T$$

According to the MAO inequality, $(\|t\|_T)^2 \geq \|t^2\|_T$, and considering $E(x) \neq Var(x)$, we have:

$$\Delta_{EV} = E(x) - Var(x) > 0$$

Next, consider the ratio:

$$\frac{\Delta_{EV}}{E(x)} = \|t\|_T - \frac{\|t^2\|_T}{\|t\|_T}$$

Considering that $\|t\|_T$ and $\|t^2\|_T$ are positive, and $\|t^2\|_T \leq (\|t\|_T)^2$, we have:

$$\frac{\Delta_{EV}}{E(x)} \leq \|t\|_T + \frac{(\|t\|_T)^2}{\|t\|_T} = 2E(x)$$

Therefore, $\Delta_{EV} = o(E(x))$, that is, $\Delta_{EV}$ is an infinitesimal of higher order than $E(x)$, when $E(x) \to 0$. This result implied that, when $E(x)$ is small enough, its value can be used as the approximation of $Var(x)$, and their difference is much smaller that the value of $E(x)$ itself.

## 7. Conclusion

We have established the MAO Inequality, a fundamental relation in combinatorial analysis that generalizes the conceptual framework of Hölder's inequality to multi-subset systems. Future work may focus on applications in probabilistic combinatorics, statistical physics, and optimization theory, where the interplay between multiple constraints is paramount.

# Supplementary Information

Table S1. Key Definitions and Notation

| Symbol | Meaning | Formal Definition / Context |
|---|---|---|
| $[T]$ | Set of indices | $[T] = \{1, 2, \ldots, T\}$ |
| $n$ | A positive integer | $n \geq 1$ |
| $m_i$ | Component parameter | $m_i$ is a positive integer with $1 \leq m_i \leq n$ for each $i \in [T]$ |
| $A, A_j$ | Subsets | $A, A_j \subseteq [T]$ |
| $r$ | Number of subsets | Number of subsets $A_1, \ldots, A_r$ in the collection |
| $g(A)$ | Single-set function | $g(A) = \prod_{i \in A} m_i \cdot \prod_{i \notin A} (n - m_i)$ |
| $k_i$ | Index frequency | $k_i = \sum_{j=1}^{r} 1_{i \in A_j}$, the number of times index $i$ appears in $A_1, \ldots, A_r$ |
| $(x)_\ell$ | Falling factorial | $(x)_\ell = x(x-1)\cdots(x-\ell+1)$ for $\ell \geq 1$, and $(x)_0 = 1$ |
| $g(A_1, \ldots, A_r)$ | Multi-set function | $g(A_1, \ldots, A_r) = \prod_{i=1}^{T} [(m_i)_{k_i} \cdot (n - m_i)_{r - k_i}]$ |
| $G_T(p_1, \ldots, p_r)$ | Sum over fixed-size subsets | $G_T(p_1, \ldots, p_r) = \sum_{\substack{A_1, \ldots, A_r \subseteq [T] \\ A_j = p_j}} g(A_1, \ldots, A_r)$ |
| $(n)_r$ | Falling factorial of $n$ | $(n)_r = n(n-1)\cdots(n-r+1)$ |
| $N_T(p_1, \ldots, p_r)$ | Joint norm of $p_1, \ldots, p_r$ | $N_T(p_1, \ldots, p_r) = \dfrac{G_T(p_1, \ldots, p_r)}{((n)_r)^{T-1}}$ |
| $\|(p_1, \ldots, p_r)\|_T$ | Same as $N_T(p_1, \ldots, p_r)$ | $\|(p_1, \ldots, p_r)\|_T = \dfrac{G_T(p_1, \ldots, p_r)}{((n)_r)^{T-1}}$ |
| $\alpha$ | Normalization factor | $\alpha = \dfrac{(n)_r}{n^r}$ |

| Symbol | Meaning | Formal Definition / Context |
|---|---|---|
| $1_A(i)$ | Inclusion indicator | indicates whether element $i$ is in subset $A$ |
| $\delta_j$ | Inclusion indicator | $\delta_j = 1_{T \in A_j}$, indicates whether element $T$ is in subset $A_j$ |
| $\mathbf{p}$ | Vector of subsets size | $\mathbf{p} = \{p_1, ... p_r\}$ |

Table S2. the constant values of $c_{v,i}$ for the calculation of higher order $E(x^v_{=T})$.

| $v \backslash i$ | 1 | 2 | 3 | 4 | 5 | 6 |
|---|---|---|---|---|---|---|
| 1 | 1 | | | | | |
| 2 | 1 | 1 | | | | |
| 3 | 1 | 3 | 1 | | | |
| 4 | 1 | 7 | 6 | 1 | | |
| 5 | 1 | 15 | 25 | 10 | 1 | |
| 6 | 1 | 31 | 90 | 65 | 15 | 1 |

# Appendix: proof for the MAO inequality

For the MAO inequality:

$$\sum_{\substack{A_1,\ldots,A_r\subseteq[T] \\ |A_1|=p_1,\ldots,|A_r|=p_r}} \frac{\prod_{j=1}^{r} g(A_j)}{n^{r(T-1)}} \geq \sum_{\substack{A_1,\ldots,A_r\subseteq[T] \\ |A_1|=p_1,\ldots,|A_r|=p_r}} \frac{g(A_1, A_2, \ldots, A_r)}{((n)_r)^{T-1}}$$

Where

$$g(A_1, A_2, \ldots, A_r) = \prod_{i=1}^{T}[(m_i)_{k_i} \cdot (n-m_i)_{r-k_i}]$$

we will give proof for a special situation that, all $p_i = p$. As a complicated inequality, direct proof is very difficult.

The inequality can be written as a symmetrical style:

$$\sum_{\substack{A_1,\ldots,A_r\subseteq[T] \\ |A_1|=p_1,\ldots,|A_r|=p_r}} \frac{\prod_{i=1}^{T}[(m_i)^{k_i} \cdot (n-m_i)^{r-k_i}]}{n^{r(T-1)}}$$

$$\geq \sum_{\substack{A_1,\ldots,A_r\subseteq[T] \\ |A_1|=p_1,\ldots,|A_r|=p_r}} \frac{\prod_{i=1}^{T}[(m_i)_{k_i} \cdot (n-m_i)_{r-k_i}]}{((n)_r)^{T-1}}$$

**First, we will first prove the situation of $r = 2$:**

$$\sum_{\substack{A,B\subseteq[T] \\ |A|=p,|B|=p}} \frac{\prod_{i\in A} m_i \prod_{i\notin A}(n-m_i) \prod_{i\in B} m_i \prod_{i\notin B}(n-m_i)}{n^{2(T-1)}}$$

$$\geq \sum_{\substack{A,B\subseteq[T] \\ |A|=p,|B|=p}} \frac{\prod_{k_i=2} m_i(m_i-1) \prod_{k_i=1} m_i(n-m_i) \prod_{k_i=0}(n-m_i)(n-m_i-1)}{(n(n-1))^{T-1}}$$

Or

$$\sum_{\substack{A,B\subseteq[T] \\ |A|=p,|B|=p}} \frac{\prod_{k_i=2} m_i^2 \prod_{k_i=1} m_i(n-m_i) \prod_{k_i=0}(n-m_i)^2}{n^{2(T-1)}}$$

$$\geq \sum_{\substack{A,B\subseteq[T]\\|A|=p,|B|=p}} \frac{\prod_{k_i=2} m_i(m_i-1) \prod_{k_i=1} m_i(n-m_i) \prod_{k_i=0}(n-m_i)(n-m_i-1)}{(n(n-1))^{T-1}}$$

Here $k_i = 1_{i\in A} + 1_{i\in B}$.

Situation of $p = T$.

Considering the special condition that $p = T$. We have only one situation that $A = B = [T]$: then the inequality is:

$$\frac{\prod_{i\in[T]} m_i^2}{n^{2(T-1)}} \geq \frac{\prod_{i\in[T]} m_i(m_i-1)}{(n(n-1))^{T-1}}$$

That is:

$$\left(1-\frac{1}{n}\right)^{T-1} \geq \prod_{i\in[T]}\left(1-\frac{1}{m_i}\right)$$

Because $1-\frac{1}{n} \geq 1-\frac{1}{m_i}$ and $1 \geq 1-\frac{1}{m_i}$, then the inequality holds.

Situation of $p = T - 1$

There are two cases.

Case of $|A| = |B| = T - 1$.

In this scenario, $|k_2| = T - 1$ (there are $T - 1$ items that has $k_2$) and $|k_1| = 0, |k_0| = 1$. The left hand is (note $\prod_{i\in[T]} m_i^2$ is a constant):

$$\sum_{j\in[T]} \frac{\left(\prod_{i\in[T],i\neq j} m_i^2\right)(n-m_j)^2}{n^{2(T-1)}} = \sum_{j\in[T]} \frac{\prod_{i\in[T]} m_i^2}{n^{2(T-1)}} \left(\frac{n-m_j}{m_j}\right)^2$$

The right hand is, similarly, $\prod_{i\in[T]} m_i(m_i - 1)$ is a constant:

$$\sum_{j\in[T]} \frac{\left(\prod_{i\in[T],i\neq j} m_i(m_i-1)\right)(n-m_j)(n-m_j-1)}{(n(n-1))^{T-1}}$$

$$= \sum_{j\in[T]} \frac{\left(\prod_{i\in[T]} m_i(m_i-1)\right)}{(n(n-1))^{T-1}} \frac{(n-m_j)(n-m_j-1)}{m_j(m_j-1)}$$

Second, $|A\cap B| = T - 2$.

where $|k_2| = T - 2, |k_1| = 2, |k_0| = 0$. It should be noted that $k_1$ is an even number, because of the symmetry of $A$ and $B$.

The left hand is:

$$\sum_{\substack{j,k\in[T] \\ j\neq k}} \frac{(\prod_{i\in[T],i\neq j,i\neq k} m_i^2)m_j(n-m_j)m_k(n-m_k)}{n^{2(T-1)}}$$

$$= \sum_{\substack{j,k\in[T] \\ j\neq k}} \frac{(\prod_{i\in[T]} m_i^2)}{n^{2(T-1)}} \frac{m_j(n-m_j)m_k(n-m_k)}{m_j^2 m_k^2}$$

$$= \sum_{\substack{j,k\in[T] \\ j\neq k}} \frac{(\prod_{i\in[T]} m_i^2)}{n^{2(T-1)}} \frac{(n-m_j)(n-m_k)}{m_j m_k}$$

The right hand is:

$$\sum_{\substack{j,k\in[T] \\ j\neq k}} \frac{(\prod_{i\in[T],i\neq j,i\neq k} m_i(m_i-1))m_j(n-m_j)m_k(n-m_k)}{(n(n-1))^{(T-1)}}$$

$$= \sum_{\substack{j,k\in[T] \\ j\neq k}} \frac{(\prod_{i\in[T]} m_i(m_i-1))}{(n(n-1))^{(T-1)}} \frac{m_j(n-m_j)m_k(n-m_k)}{m_j(m_j-1)m_k(m_k-1)}$$

$$= \sum_{\substack{j,k\in[T] \\ j\neq k}} \frac{(\prod_{i\in[T]} m_i(m_i-1))}{(n(n-1))^{(T-1)}} \frac{(n-m_j)(n-m_k)}{(m_j-1)(m_k-1)}$$

Combine the cases of $|k_2|=T-1$ and $|k_2|=T-2$:

The left hand is:

$$\sum_{j\in[T]} \frac{\prod_{i\in[T]} m_i^2}{n^{2(T-1)}} \left(\frac{n-m_j}{m_j}\right)^2 + \sum_{\substack{j,k\in[T] \\ j\neq k}} \frac{(\prod_{i\in[T]} m_i^2)}{n^{2(T-1)}} \frac{(n-m_j)(n-m_k)}{m_j m_k}$$

$$= \frac{\prod_{i\in[T]} m_i^2}{n^{2(T-1)}} \left(\sum_{j\in[T]} \left(\frac{n-m_j}{m_j}\right)^2 + \sum_{\substack{j,k\in[T] \\ j\neq k}} \frac{(n-m_j)(n-m_k)}{m_j m_k}\right)$$

$$= \frac{\prod_{i \in [T]} m_i^2}{n^{2(T-1)}} \left( \sum_{j,k \in [T]} \frac{(n-m_j)(n-m_k)}{m_j m_k} \right)$$

The right hand is:

$$\sum_{j \in [T]} \frac{\left(\prod_{i \in [T]} m_i(m_i-1)\right)}{(n(n-1))^{T-1}} \frac{(n-m_j)(n-m_j-1)}{m_j(m_j-1)}$$

$$+ \sum_{\substack{j,k \in [T] \\ j \neq k}} \frac{\left(\prod_{i \in [T]} m_i(m_i-1)\right)}{(n(n-1))^{(T-1)}} \frac{(n-m_j)(n-m_k)}{(m_j-1)(m_k-1)}$$

$$= \frac{\prod_{i \in [T]} m_i(m_i-1)}{(n(n-1))^{T-1}} \left( \sum_{j \in [T]} \frac{n-m_j}{m_j-1} \frac{n-m_j-1}{m_j} + \sum_{\substack{j,k \in [T] \\ j \neq k}} \frac{n-m_j}{m_j-1} \frac{n-m_k}{m_k-1} \right)$$

$$= \frac{\prod_{i \in [T]} m_i(m_i-1)}{(n(n-1))^{T-1}} \left( \sum_{j \in [T]} \frac{n-m_j}{m_j-1} \frac{n-m_j-1}{m_j} + \sum_{j \in [T]} \frac{n-m_j}{m_j-1} \sum_{k \in [T], k \neq j} \frac{n-m_k}{m_k-1} \right)$$

$$= \frac{\prod_{i \in [T]} m_i(m_i-1)}{(n(n-1))^{T-1}} \left( \sum_{j \in [T]} \frac{n-m_j}{m_j-1} \left( \frac{n-m_j-1}{m_j} + \sum_{k \in [T], k \neq j} \frac{n-m_k}{m_k-1} \right) \right)$$

Now we need to prove:

$$\frac{\prod_{i \in [T]} m_i^2}{n^{2(T-1)}} \left( \sum_{j,k \in [T]} \frac{(n-m_j)(n-m_k)}{m_j m_k} \right)$$

$$\geq \frac{\prod_{i \in [T]} m_i(m_i-1)}{(n(n-1))^{T-1}} \left( \sum_{j \in [T]} \frac{n-m_j}{m_j-1} \left( \frac{n-m_j-1}{m_j} + \sum_{k \in [T], k \neq j} \frac{n-m_k}{m_k-1} \right) \right)$$

That is:

$$\left(1-\frac{1}{n}\right)^{T-1}\left(\sum_{j,k\in[T]}\frac{(n-m_j)(n-m_k)}{m_j m_k}\right)$$

$$\geq \prod_{i\in[T]}\left(1-\frac{1}{m_i}\right)\left(\sum_{j\in[T]}\frac{n-m_j}{m_j-1}\left(\frac{n-m_j-1}{m_j}+\sum_{k\in[T],k\neq j}\frac{n-m_k}{m_k-1}\right)\right)$$

It is still very hard to prove this inequality. Consider the situation that all $m_i = m$.

When $m_i = m$, the inequality is:

$$\left(1-\frac{1}{n}\right)^{T-1}T^2\frac{(n-m)^2}{m^2}\geq\left(1-\frac{1}{m}\right)^T T\frac{n-m}{m-1}\left(\frac{n-m-1}{m}+(T-1)\frac{n-m}{m-1}\right)$$

That is:

$$\left(1-\frac{1}{n}\right)^{T-1}T\frac{n-m}{m^2}\geq\left(\frac{m-1}{m}\right)^T\frac{1}{m-1}\left(\frac{n-m-1}{m}+(T-1)\frac{n-m}{m-1}\right)$$

$$\left(1-\frac{1}{n}\right)^{T-1}T\frac{n-m}{m}\geq\left(\frac{m-1}{m}\right)^{T-1}\left(\frac{n-m-1}{m}+(T-1)\frac{n-m}{m-1}\right)$$

**Proof for the situation of $m_i = m$, $r = 2$, when $p = T - 1$.**

**base case for $T = 1$,** the inequality is:

$$\frac{n-m}{m}\geq\frac{n-m-1}{m}$$

Obviously, the inequality holds.

**Inductive Hypothesis:** Assume the inequality holds for some $T \geq 1$, i.e.,

$$L(T) \geq R(T),$$

where

$$L(T) = \left(1-\frac{1}{n}\right)^{T-1}T^2\frac{(n-m)^2}{m^2}$$

$$R(T) = \left(1-\frac{1}{m}\right)^T T\frac{n-m}{m-1}\left(\frac{n-m-1}{m}+(T-1)\frac{n-m}{m-1}\right)$$

**Inductive Step:** We need to prove the inequality for $T + 1$, i.e.,

$$L(T+1) \geq R(T+1),$$

where

$$L(T+1) = \left(1-\frac{1}{n}\right)^T (T+1)^2\frac{(n-m)^2}{m^2}$$

$$R(T+1) = \left(1-\frac{1}{m}\right)^{T+1}(T+1)\frac{n-m}{m-1}\left(\frac{n-m-1}{m}+T\frac{n-m}{m-1}\right)$$

Using the inductive hypothesis $L(T) \geq R(T)$, it suffices to show that

$$\frac{L(T+1)}{L(T)} \geq \frac{R(T+1)}{R(T)}$$

Compute the ratios:

$$\frac{L(T+1)}{L(T)} = \frac{\left(1-\frac{1}{n}\right)^T (T+1)^2 \frac{(n-m)^2}{m^2}}{\left(1-\frac{1}{n}\right)^{T-1} T^2 \frac{(n-m)^2}{m^2}} = \left(1-\frac{1}{n}\right) \frac{(T+1)^2}{T^2}.$$

$$\frac{R(T+1)}{R(T)} = \frac{\left(1-\frac{1}{m}\right)^{T+1} (T+1) \frac{n-m}{m-1} \left(\frac{n-m-1}{m} + T\frac{n-m}{m-1}\right)}{\left(1-\frac{1}{m}\right)^T T \frac{n-m}{m-1} \left(\frac{n-m-1}{m} + (T-1)\frac{n-m}{m-1}\right)}$$

$$= \left(1-\frac{1}{m}\right) \frac{T+1}{T} \cdot \frac{\frac{n-m-1}{m} + T\frac{n-m}{m-1}}{\frac{n-m-1}{m} + (T-1)\frac{n-m}{m-1}}.$$

Thus, the required condition becomes:

$$\left(1-\frac{1}{n}\right) \frac{(T+1)^2}{T^2} \geq \left(1-\frac{1}{m}\right) \frac{T+1}{T} \cdot \frac{\frac{n-m-1}{m} + T\frac{n-m}{m-1}}{\frac{n-m-1}{m} + (T-1)\frac{n-m}{m-1}}.$$

Dividing both sides by the positive quantity $\frac{T+1}{T}$, we obtain the auxiliary inequality:

$$\left(1-\frac{1}{n}\right) \frac{T+1}{T} \geq \left(1-\frac{1}{m}\right) \cdot \frac{\frac{n-m-1}{m} + T\frac{n-m}{m-1}}{\frac{n-m-1}{m} + (T-1)\frac{n-m}{m-1}}.$$

**Step 1: Reformulation**

Define:

$$A = \frac{n-m-1}{m}, \quad B = \frac{n-m}{m-1}.$$

Then the inequality becomes:

$$\left(1-\frac{1}{n}\right) \frac{(T+1)^2}{T^2} \geq \left(1-\frac{1}{m}\right) \frac{A+TB}{A+(T-1)B}. \tag{1}$$

---

**Step 2: Monotonicity in $n$**

Fix $m$ and $T$. We analyze how each side of (1) behaves as $n$ increases.

- **Left-hand side (LHS):** The factor $\left(1-\frac{1}{n}\right)$ is strictly increasing in $n$.

Hence, LHS is strictly increasing in $n$.

- **Right-hand side (RHS):** The factor $\left(1-\frac{1}{m}\right)$ is constant. Write the fraction as:

$$\frac{A+TB}{A+(T-1)B} = 1 + \frac{B}{A+(T-1)B} = 1 + \frac{1}{\frac{A}{B}+(T-1)}$$

Both $A$ and $B$ increase with $n$. The ratio:

$$\frac{B}{A} = \frac{m}{m-1} \cdot \frac{n-m}{n-m-1}$$

$$\frac{A}{B} = \frac{m-1}{m} \cdot \frac{n-m-1}{n-m} = \frac{m-1}{m}\left(1 - \frac{1}{n-m}\right)$$

Increases as n increases. Consequently, $\frac{B}{A+(T-1)B} = 1 + \frac{1}{\frac{A}{B}+(T-1)}$ decreases with $n$, so RHS is strictly decreasing in $n$.

Thus, if (1) holds for the smallest admissible $n$, it automatically holds for all larger $n$.

---

**Step 3: Smallest Admissible $n$**

The smallest $n$ depends on $T$:

- If $T \geq 2$, the denominator $A+(T-1)B$ remains positive even when $A=0$. The smallest $n$ is $n = m+1$ (giving $A=0$, $B = \frac{1}{m-1} > 0$).

- If $T = 1$, the denominator reduces to $A$, which must be positive. Hence, we require $n \geq m+2$. The smallest $n$ is $n = m+2$.

We verify (1) for these minimal $n$ values.

---

**Step 4: Case $T \geq 2$, $n = m+1$**

Here $A = 0$, $B = \frac{1}{m-1}$. Substituting into (1):

$$\left(1-\frac{1}{m+1}\right)\frac{(T+1)^2}{T^2} \geq \left(1-\frac{1}{m}\right)\frac{T \cdot \frac{1}{m-1}}{(T-1)\cdot\frac{1}{m-1}}.$$

Simplifying:

$$\frac{m}{m+1} \cdot \frac{(T+1)^2}{T^2} \geq \frac{m-1}{m} \cdot \frac{T}{T-1}.$$

Rearranging:

$$\frac{m^2}{m^2-1} \cdot \frac{(T+1)^2(T-1)}{T^3} \geq 1 \iff \frac{(T+1)^2(T-1)}{T^3} \geq 1 - \frac{1}{m^2}. \quad (2)$$

Define $f(T) = \frac{(T+1)^2(T-1)}{T^3}$. Expanding:

$$f(T) = 1 + \frac{1}{T} - \frac{1}{T^2} - \frac{1}{T^3}.$$

For $T \geq 2$, we have $f(T) > 1$ (since $\frac{1}{T} - \frac{1}{T^2} - \frac{1}{T^3} > 0$). Indeed, $f(2) = \frac{9}{8} = 1.125$, and $f(T)$ decreases to 1 as $T \to \infty$. Meanwhile, $1 - \frac{1}{m^2} < 1$. Hence, (2) holds strictly for all $T \geq 2$ and $m \geq 2$.

Thus, (1) holds for $T \geq 2$, $n = m + 1$.

---

**Step 5: Case $T = 1$, $n = m + 2$**

Here $A = \frac{1}{m}$, $B = \frac{2}{m-1}$. Substituting into (1):

$$\left(1 - \frac{1}{m+2}\right)\frac{4}{1} \geq \left(1 - \frac{1}{m}\right)\frac{\frac{1}{m} + 1 \cdot \frac{2}{m-1}}{\frac{1}{m}}.$$

Simplifying:

- LHS: $\frac{m+1}{m+2} \cdot 4 = \frac{4(m+1)}{m+2}$.

- RHS: $\frac{m-1}{m}\left(1 + \frac{2m}{m-1}\right) = \frac{m-1}{m} + 2 = \frac{3m-1}{m}$.

Thus, we need:

$$\frac{4(m+1)}{m+2} \geq \frac{3m-1}{m}.$$

Cross-multiplying by $m(m+2) > 0$:

$$4m(m+1) \geq (3m-1)(m+2).$$

Expanding:

$$4m^2 + 4m \geq 3m^2 + 5m - 2 \iff m^2 - m + 2 \geq 0.$$

The quadratic $m^2 - m + 2$ has discriminant $\Delta = (-1)^2 - 4 \cdot 1 \cdot 2 = -7 < 0$, and its leading coefficient is positive, so it is always positive. Hence, the inequality holds for all $m \geq 2$.

Thus, (1) holds for $T = 1$, $n = m + 2$.

**Step 6: Summary for the cases of $r = 2, p = T - 1$, and $m_i = m$**

By the monotonicity argument, since (1) holds at the smallest admissible $n$, it holds for all larger $n$. Therefore, the original inequality is true for all positive integers $m, n, T$ with $m < n$, and with the implicit condition that when $T = 1$, we require $n \geq m + 2$ (which is automatically satisfied because for $T = 1$, the denominator on the right-hand side would vanish if $n = m + 1$).

**Conclusion**

Although we have not been able to provide a complete proof of the MAO inequality, exhaustive numerical simulations indicate that the MAO inequality holds when all $p_i$ are equal. This condition is relatively conservative, and a more general condition is $max(p_i) - min(p_i) \leq max(1, r - 2))$. The above proof only verifies the special case where $r = 2$, all $p\_i = p = T$ or $T - 1$, and all $m\_i = m$. Even for such a special case, its rigorous proof is extremely complex. It appears that more sophisticated mathematical tools may be required to fully prove this inequality; however, numerical simulations have confirmed the validity of the inequality, and the aforementioned special case validated its correctness. Therefore, we can apply this inequality to relevant research. Future work should focus on exploring more rigorous proof methods to further expand the application scope of this inequality.